\documentclass[useAMS]{gCOM2e}

\usepackage{url}
\usepackage{algorithm}
\usepackage{algorithmic}
\usepackage{fancyvrb}

\begin{document}


\articletype{Research Paper}

\title{A PCG Implementation of an Elliptic Kernel in an Ocean Global Circulation Model Based on GPU Libraries}


\author{S. Cuomo$^{\rm a}$$^{\ast}$, P. De Michele$^{\rm a}$ and R. Farina$^{\rm b}$\thanks{$^\ast$Corresponding author. Email: salvatore.cuomo@unina.it\vspace{6pt}} and M. Chinnici$^{\rm c}$\\\vspace{6pt}  $^{\rm a}${\em{University of Naples ``Federico~II'' - Dept. of Mathematics and Applications ``R.~Caccioppoli'' - Via Cinthia, 80126, Naples, Italy}}; $^{\rm b}${\em{CMCC, Research Center, Via Aldo Moro 44 (Bologna), Italy}} $^{\rm c}${\em{ENEA-UTICT, Casaccia Research Center, S.~Maria di Galeria (Roma), Italy}}\\\vspace{6pt}\received{released October 2012} }

\maketitle

\begin{abstract}In this paper an inverse preconditioner for the numerical solution of an elliptic Laplace problem of a global circulation ocean model is presented. The inverse preconditiong technique is adopted in order to efficiently compute the numerical solution of the elliptic kernel by using the Conjugate Gradient (CG) method. We show how the performance and the rate of convergence of the solver are linked to the discretized grid resolution and to the Laplace coefficients of the oceanic model. Finally, we describe an easy-to-implement version of the solver on the Graphical Processing Units (GPUs) by means of scientific computing libraries and we discuss its performance.\bigskip

\begin{classcode}65Y5; 65Y10; 65F08; 65F35; 37N10.\end{classcode}\bigskip
\end{abstract}

\section{Introduction}
The ocean is now well known to play a dominant role in the climate system because it can initiate and amplify climate change on many different time scales. Hence, the simulation of ocean model is became a relevant but highly complex task and it involves an intricate interaction of theoretical insight, data handling and numerical modelling. 
Over the past several years, ocean numerical models have become quite realistic as a result of improved methods, faster computers, and global data sets. Models now treat basin-scale to global domains while retaining the fine spatial scales that are important for modelling the transport of heat, salt, and other properties over vast distances.
Currently, there are many models and methods employed in the rapidly advancing field of numerical ocean circulation modelling as Nemo, Hops, MOM , POP et al (see~\cite{paper1} for a nice review). However, several of these numerical models are not yet optimized by using scientific computing libraries and ''ad hoc'' preconditioning techniques. In all these frameworks the numerical kernel is represented  by the discretization of the Navier-Stokes equations ~\cite{paper3} on a three dimensional grid and by the computation of the evolution time of each variable for each grid point. The high resolution computational grid requires efficient preconditiong techniques for improving the accuracy in the computed solution and  parallelization strategies for answering to the huge amount of computational demand.\\
In this paper we propose a new solver based on preconditioned conjugate gradient (PCG) method with an approximate inverse preconditioner  AINV~\cite{book4-14}   for the numerical solution of the elliptic sea-surface equation in NEMO-OPA ocean model~\cite{book3-1}, a state of the art modelling framework in the oceanographic research. The PCG is  a widely used iterative method for solving linear systems with symmetric, positive definite matrix and it has proven its efficiency and robustness in a lot of applications. The preconditioning is often a bottleneck in solving the linear systems efficiently and it is well established that a suitable preconditioner increases the performance of an application dramatically.\\
The elliptic sea-surface equation  is originally solved in NEMO-OPA by means of the PCG with diagonal preconditioner, and  in our work we prove  that it is  inefficient and inaccurate. We build a new inverse preconditioner and we implement the PCG on a Graphic Processor Unit (GPU) by means of the linear algebra Scientific Computing libraries. The GPUs are massively parallel architectures that efficiently work with the linear algebra operations and give impressive performance improvements. They require a deep understanding of the underlying computing architecture and the programming with these devices involves a massive re-thinking of existing CPU based applications. A challenge is how to optimally use the GPU hardware adopting adequate programming techniques, models, languages and  tools. In this paper, we present an easy-to-implement version of the elliptic solver with the scientific computing libraries on Compute Unified Device Architecture (CUDA)~\cite{book4-2bis2}. We implement a code  by using  CUDA  based supported libraries  CUBLAS~\cite{book4-2} and CUSPARSE~\cite{book4-2bis} for the sparse linear algebra and the graph computations.
\noindent The library GPU based approach allows a short code development times and an easy to use  GPUs implementations that  can fruitfully  speed up  the expensive numerical kernel of an oceanographic simulator. The paper is set out as follows. In section 2 we briefly review the mathematical model: elliptic equations that are at the heart of the model. In section 3, the preconditioned conjugate gradient method used to invert the elliptic equations are described. In section 4, we outline a implementation strategy for solving the elliptic solver by using standard libraries and in section 5 we discuss the mapping of our algorithm onto a massively parallel machine. Finally, the conclusions are drawn.  

\section{The Mathematical Model}
Building and running ocean models able to simulate the world of global circulation  with great realism require a variety of scientific skills. In modelling the general ocean circulation it is necessary to solve problems of elliptic nature. These problems are difficult to solve, with the following issues causing the most trouble in practice\cite{paper1}:
\begin{enumerate}
\item {In simulations with complicated geometry (e.g., multiple islands), topography, time varying
surface forcing, and many space-time scales of variability (i.e., the World Ocean), achieving a
good first guess for the iterative elliptic solver is often quite difficult to achieve. This makes it
difficult for elliptic solvers to converge to a solution within a reasonable number of iterations. For this reason, many climate modellers limit the number of elliptic solver iterations used, even if the solver has not converged. This approach is very unsatisfying.}

\item{Many elliptic solvers with their associated non-local and time dependent boundary conditions (be they Neumann or Dirichlet) do not project well onto parallel distributed computers, which acts to hinder their scaling properties \cite{duco,web,webi}.}
\end{enumerate}

\begin{figure}[h!]
	\centering
	\includegraphics[scale=0.25]{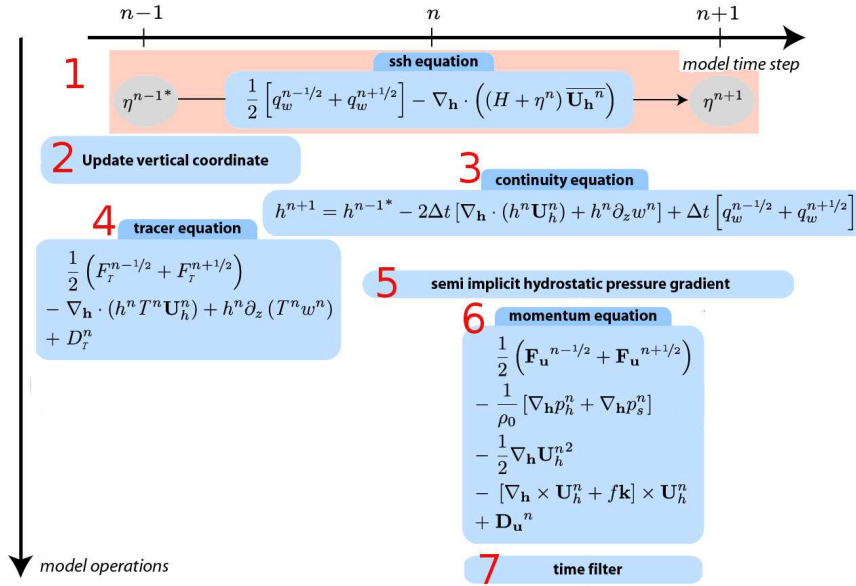}
	\caption{\small {NEMO-OPA model.}}
	\label{fig:NEMO-OPA}
\end{figure}

\noindent In the OPA-NEMO numerical code the primitive equations are discretized within sea-surface hypothesis~\cite{Roullet} and  the model is charecterized by the three-dimensional distribution of currents, potential temperature, salinity, preassure and density~\cite{paper3}. 
The numerical method OPA-NEMO is grounded on discretizing of the primitive equation - by the use of  finite differences on a three dimensional grid - and computing the time evolution to each variable "ocean"  at each grid point for the entire globe~\cite{ara}.
A sketch of the  OPA-NEMO computational model, see Figure~\ref{fig:NEMO-OPA}, shows the complex dynamic processes that mimic the ocean circulation model, composed by steps that are many time simulated.

\noindent  The kernel algorithm (highlighted with red color, equation 1) solves the sea-surface hight equation $\eta$ The elliptic kernel is discretized  with a semi-discrete equations, as following:
\begin{center}
	\begin{small}
		\begin{eqnarray}
			{\eta^{n+1}=\eta^{n-1}}- 2 \Delta t D^n  \qquad \qquad \qquad \quad \qquad \ \ \label{uno} \quad \\[1.5mm]
			2\Delta t g T_{c} \Delta_h D^{n+1}= D^{n+1}- D^{n-1}+2\Delta t g \Delta_h \eta^n \label{unounos} \quad \ \ \\[1.5mm]
			\Delta_h = \nabla \big[(H+\eta^n)\nabla \big]. \qquad \qquad \quad \qquad \ \ \qquad \quad
		\end{eqnarray}
	\end{small}
\end{center}

\noindent  where $\eta^n, \ n  \in \mathbb{N}$ is  the sea-surface height at the $n-$th step of the model,  which describes the shape of the air-sea interface, $D^n$ is the  centered difference approximation of the first time derivative of $\eta$,
   $\Delta t $ is the time stepping, $g$ is the gravity constant, $T_c$ is a physical parameter, $\Delta_h$ is the horizontal Laplacian operator and  $H$  is the depth
of the ocean bottom  \cite{book3-1}. 
\noindent Whereas the domain of the ocean models is the Earth sphere (or part of it) the model uses  the geographical coordinates system $(\lambda,\phi,r)$ in which a position is defined by the latitude $\phi$, the longitude $\lambda$ and the distance from the center of the earth $r=a + z(k)$ where a is the Earth's radius and z the altitude above a reference sea level. The local deformation of the curvilinear geographical  coordinate system is given by  $e_1$,$e_2$ and $e_3$:

\begin{equation}
	\begin{small}
		\begin{array}{c}
		e_1  =  rcos \phi, \quad 
		e_2  =  r, \quad 
		e_3  = 1.
		\end{array}
	\end{small}
\end{equation}

The Laplacian  Operator in spherical coordinates  $\Delta_h D^{n+1}$ in~(\ref{unounos})  becomes:

\begin{small}
{
	\begin{equation}
		\Delta_h D^{n+1}=\frac{1}{e_1 e_2} \Bigg[\frac {\partial }{\partial i  } \bigg( \alpha(\phi )\frac {\partial D^{n+1} }{\partial i  } \bigg)+\frac {\partial }{\partial j  } \bigg(\beta(\phi)\frac {\partial D^{n+1} }{\partial i  } \bigg) \Bigg]
	\end{equation}}
\end{small}

\noindent where:

\begin{eqnarray} 
\alpha(\phi)= (H+\eta^n) {{  e_2}/{ e_1}} \label{uu} \\[2mm]
\beta(\phi)= (H+\eta^n)  {e_1}/{ e_2}  \label{ua} 
\end{eqnarray}

\noindent For  the functions  $\alpha(\phi)$ in (\ref{uu}) and $\beta(\phi)$ in (\ref{ua}), we have the following relations:

\begin{equation}
	\begin{small}
		\lim_{\phi \longrightarrow \pm \frac{\pi}{2} }\alpha(\phi)=+\infty \quad \wedge \
		\lim_{\phi \longrightarrow \pm \frac{\pi}{2} }\beta(\phi)=0
		\label{limite1}
	\end{small}
\end{equation}

\noindent From~(\ref{limite1}),  if we choose  $M, \epsilon \in \mathbb{R}$ with  $M>>\epsilon$ then exists an interval $\big[\frac{\pi}{2}-\delta,  \frac {\pi}{2}\big]$ or $\big[-\frac{\pi}{2}, \ -\frac{\pi}{2}+~\delta\big]$, such that the following inequality holds:

\begin{equation}
	\begin{small}
		\alpha(\phi)>M>>\epsilon>\beta(\phi)
		\label{limite2}
	\end{small}
\end{equation}

\noindent In physical terms, in the proximity of the geographical poles, $(\lambda, \pm \phi/2,r)$, there are several orders of magnitude between the functions $\alpha(\phi)$ and $\beta(\phi)$

\noindent The result  (\ref{limite2}), will  significantly influence the rate of convergence in the iterative solver.
	
\section{Inverse Preconditioned Techniques in the Elliptic Solver of the Ocean Model}

\noindent Let us now consider the elliptic NEMO model~\cite{book3-1} defined by the following coefficients:
\begin{equation}
	\begin{small}
		\begin{array}{c}
			C_{i,j}^{NS}= 2\Delta t^2 { H(i,j) e_{1}(i,j)}/{e_{2}(i,j)},  \quad 
			C_{i,j}^{EW}= 2\Delta t^2  { H(i,j) e_{2}(i,j)}/{e_{1}(i,j)}  \\[0.5mm]
			b_{i,j}=\delta_i(e_{2}M_u)-\delta_j(e_{1}M_v) \quad \quad 
		\end{array}
		\label{coef}
	\end{small}
\end{equation}
where $\delta_i $ and $\delta_j$ are the discrete derivative operators along the  axes $\mathbf i$ and $\mathbf j$.
\noindent The discretization of the equation~(\ref{unounos}) by means of a five-point finite difference method gives:
\begin{equation}
	\begin{small}
		\begin{array}{c}
			C_{i,j}^{NS} D_{i-1,j}+C_{i,j}^{EW} D_{i,j-1}-
			\big(C_{i+1,j}^{NS}+C_{i,j+1}^{EW}+ C_{i,j}^{NS}+ C_{i,j}^{EW}\big) D_{i,j}+\\[2mm]
			\qquad \qquad + C_{i,j+1}^{EW} D_{i,j+1} +  C_{i+1,j}^{NS} D_{i+1,j}=b_{i,j}. 
		\end{array}
		\label{sistema}
	\end{small}
\end{equation}

\noindent where the equation~(\ref{sistema}) is a symmetric system of linear equations. All the elements of the sparse matrix $\mathbf A$ vanish except those of five diagonals. With the natural ordering of the grid points (i.e. from west to east and from south to north), the structure of $\mathbf A$ is a block-tridiagonal with tridiagonal or diagonal blocks. The matrix $\mathbf A$ is a positive-definite symmetric matrix with $n=jpi \times jpj$ size, where $jpi$ and $jpj$ are respectively the horizontal dimensions of the grid discretization of the domain.\\ The Conjugate Gradient Method  is a very efficient iterative method for solving the linear system~(\ref{sistema}) and it  provides the exact solution in a number of iterations equal to the size of the matrix. The convergence rate is faster as the matrix is closer to the identity one. By spectral point of view a convergence relation between the solution of the linear system and its approximation $x_{m}$ is given by:

\begin{small}
	\begin{equation}
		\|\mathbf x-\mathbf x_{m}\|_A< 2\bigg( \frac {\sqrt {\mu_2(A)} -1}{\sqrt {\mu_2(A)} + 1}\bigg)^{m-1}\|\mathbf x-\mathbf x_{0}\|_A
		\label{rel}
	\end{equation}
\end{small}

\noindent with $\mu_2(A)=\lambda_{max}/\lambda_{min}$, where $\lambda_{max}$ and $\lambda_{min}$ are respectively the greatest and the lowest eigenvalue of $\mathbf A$, and $\|\cdot\|_A$ is the A-norm. The preconditioning framework consists to introduce amatrix $ \mathbf M$, that is an approximation of $\mathbf A $ easier to invert, and to solve the equivalent linear system:
\begin{equation}
	\begin{small}\mathbf{M^{-1}A}\mathbf x=\mathbf {M^{-1}b}\end{small}
\end{equation}
The ocean global model NEMO-OPA uses the diagonal preconditioner, where $\mathbf M$ is chosen to the diagonal of $\mathbf A$. Let us introduce the following cardinal coefficients:

\begin{small}
\begin{eqnarray}
			\alpha_{i,j}^E={{C_{i,j+1}}^{EW}}/{d_{i,j}} \qquad \alpha_{i,j}^W= {{C_{i,j}}^{EW}}/{d_{i,j}}   \qquad  \label{coef21} \\[2mm] \beta_{i,j}^S= {{C_{i,j}}^{NS}}/{d_{i,j}} \qquad \beta_{i,j}^N= {{C_{i+1,j}}^{NS}}/{d_{i,j}}
		\label{coef22}
	\end{eqnarray}
\end{small}
\noindent where $d_{i,j}= \big(C_{i+1,j}^{NS}+C_{i,j+1}^{EW}+ C_{i,j}^{NS}+ C_{i,j}^{EW}\big)$ represents the diagonal of the matrix $\mathbf A$. The~(\ref{sistema}), using the diagonal preconditioner, can be written as:
\begin{equation}
	\begin{small}
		\begin{array}{c}
			-\beta_{i,j}^{S} D_{i-1,j}-\alpha_{i,j}^{W} D_{i,j-1}+ D_{i,j}
			- \alpha_{i,j}^{E} D_{i,j+1}+
			\hfill - \beta_{i,j}^{N} D_{i+1,j} =\bar b_{i,j}. 
		\end{array}
		\label{nuovosistema}
	\end{small}
\end{equation}
\noindent with $\bar b_{i,j}= - b_{i,j} /d_{i,j}$.

\noindent Starting from the observations~(\ref{limite1}) and~(\ref{limite2}) we proof that the diagonal preconditioner does not work very well in some critical physical situations involving curvilinear spherical coordinates.

\begin{proposition}
\label{remark}
 In the geographical coordinate, if $\phi \rightarrow + {\frac {\pi}{2}}^{-}$, $\Delta \lambda \rightarrow 0 $,  $\Delta \phi \rightarrow 0$ then the conditioning number $\mu (\mathbf {M^{-1}A)}$ goes to $+\infty$. 

\end{proposition}

\begin{small}

\vspace{0.5cm}

\noindent \textbf{Proof.} \ \   In the geographical coordinate, i.e. when $(i, j) \rightarrow (\lambda, \phi)$ and $(e_1 , e_2 ) \rightarrow (r \cos \phi, r)$,  for $\phi \rightarrow + {\frac {\pi}{2}}^{-}$, $\Delta \lambda \rightarrow 0 $ and $\Delta \phi \rightarrow 0$, the functions  $\alpha^{W}$ and $\alpha^{E}$ in (\ref{coef21})  go to -1/2 while $\beta^{N}$ and $\beta^{S}$ in (\ref{coef22}) go to 0. Hence the limit of matrix $\mathbf {M^{-1}A}$ is given by:

\begin{equation}   
\begin{tiny}
\mathbf A'=\left [
\begin{array}{ccccc}
1& - 1/2 & 0& \ldots & 0  \\
 -1/2 & 1 & -1/2 & \ddots& 0\\
0 & \ddots & \ddots & \ddots & 0\\
\vdots & \ddots & \ddots & \ddots &-1/2\\
0 & \ddots & \ddots & -1/2 & 1\\
\end{array}
\right ].
\end{tiny} 
\label{matricelimite}
\end{equation}

\vspace{0.5cm}

\noindent The eigenvalues of the matrix in (\ref{matricelimite}) are:

\begin{equation}
\lambda_k=1+\cos\bigg(\frac {k \pi}{n+1}\bigg) \quad k=1,...,n
\end{equation}
and then  the condition number $\mu_2(M^{-1}{A}) =\lambda_{max}/\lambda_{min} \approx n^2/2$ (by using the series expansion of $ \cos x = 1-x ^ 2/2 + o (x ^ 2) $). Moreover for  $\Delta \lambda \rightarrow 0 $ and $\Delta \phi \rightarrow 0$ the size $n$ of the matrix $\mathbf A $  goes to $ +\infty$ and hence we obtain the thesis. { \hfill  $ \small \blacksquare $ }
\end{small} \\[2mm]

\noindent By the proposition~(\ref{remark}), for $n$ large and $\phi \rightarrow \pm \pi/2$, it is preferable to adopt more suitable preconditioning techniques or a stategy based on the local change of the coordinates at poles. In this paper we propose an alternative approximate sparse inverse preconditioning AINV techniques \cite{book4-14} for the linear system (\ref{sistema}).
AINV technique is a critical step since the inverse of a sparse matrix is usually
dense. The problem is how to build a preconditoner that preserves the sparse structure. We introduce  a factored sparse approximate inverse FSAI preconditioner  $\mathbf P= \tilde {\mathbf Z}  \tilde {\mathbf {Z^t}}$ \cite{book4-16,book4-31}, computed by means conjugate-orthogonalization procedure.  Specifically, we propose  an ``ad hoc'' method for computing an incomplete factorization of the inverse of the matrix $\mathbf T \subset \mathbf A$,  obtained by $\mathbf A$  taking only the elements $a_{i,j}$  such that $|j-i|\leq 1$. The factorized sparse approximate inverse of $\mathbf T$ is used as explicit preconditioner  for (\ref{sistema}).
\noindent In the following we give some remarks for the sparsity pattern selection S of our inverse preconditioner $\mathbf P$.

\begin{proposition}{}
\label{titi}
If $\mathbf T$ is a tridiagonal, symmetric and diagonally dominant matrix, with diagonal elements all positive $t_{k,k}>0,\ k=1,...,n$, then the Cholesky's factor $\mathbf U$ of the matrix  $\mathbf T$  is again diagonally dominant.
\end{proposition}

\begin{small}
\vspace{0.2cm}
\noindent \textbf{Proof.} Since $\mathbf T$  is a tridiagonal matrix then  $\mathbf  U$ is a bidiagonal matrix. 
Using the inductive method we  proof  that $\mathbf U$ is diagonally dominant matrix. For $k=1$ is trivially, indeed  by hypothesis we know that  $|a_{1,1}|>|a_{1,2}| \Longleftrightarrow  |u_{1,1}^2|> |u_{1,1} u_{1,2}|$, then we  obtain  $|u_{1,1}|>|u_{1,2}|$. Moreover placed the thesis true for $k-1$ i.e. $|{{u_{k-1,k-1}}}|> |{{u_{k-1,k}}}|$ then by the following inequalities:
\begin{eqnarray}
|a_{k,k}|>|a_{k,k-1}|+|a_{k,k+1}|  \Longleftrightarrow  \qquad \quad  \qquad \quad \qquad \quad \qquad \quad  \nonumber \\[2mm] |u_{k-1,k}^2+ u_{k,k}^2|>|u_{k-1,k} u_{k-1,k-1} |+|u_{k,k} u_{k,k+1}| 
 > u_{k-1,k}^2+|u_{k,k} u_{k,k+1}|.  \qquad \quad  \label{dividi}
 \end{eqnarray}
   subtracting  the inequality (\ref{dividi}) for  $u_{k-1,k}^2$, then the thesis holds also for k.{ \hfill  $ \small \blacksquare $ }
\end{small}

\vspace{2mm}

\noindent This result allows to prove the following proposition:

\vspace{2mm}

\begin{proposition}
The inverse matrix $\mathbf Z$ of a bidiagonal and diagonally dominant matrix  $\mathbf U$ has column vectors $\mathbf z_k, k=1,...n$ such that starting from diagonal element $z_{k,k}$, they  contain a  finite sequence $ \{z_{k-i,k} \}_{i=0,...,k-1}$  strictly decreasing.
\label{final}
\end{proposition} 

\begin{small}
\vspace{0.2cm}
\noindent \textbf{Proof.}   Applying a  backward substitution procedure for solving the system of equations $\mathbf U \mathbf z_k = \mathbf e_k $,  we get:

\begin{equation}
z_{k-i,k}=\left\{
\begin{array}{l}
{1}/{u_{k,k}} \quad  {if \ \ i=0}\\[5mm]
({-1})^{i}/{u_{k,k}} \cdot\ {\displaystyle  \prod_{r=1}^{i}}
\big({u_{k-r,k-r+1}}/{u_{k-r,k-r}}\big)\\[5mm]
 \qquad \qquad \qquad \qquad \qquad { if \ \  0<i\leq k-1}.
\end{array}\right.
\end{equation}
\label{propdue}  
By means of the preposition  (\ref{titi}) we obtain that $z_{k-i,k}> z_{k-i-1,k}$  with  {$ \small{0<i\leq k-1}$} and hence the thesis is proved. { \hfill  $ \small \blacksquare $ }

\end{small}

\vspace{2mm}





\noindent The previous propositions (\ref{titi}) and (\ref{final}) enable  to select a sparsity  pattern S by the following scheme: 

\begin{enumerate}
{\tt \small
\item{ Consider the symmetric, diagonally dominant and triangular matrix $\mathbf{T}$,  obtained by  $\mathbf A $  taking only the elements $a_{i,j}$ such that $|j-i|\leq 1$}
\item{ $\mathbf T= \mathbf U ^T \mathbf U$ is diagonally dominant matrix. Consequently its Cholesky factor  $\mathbf U$ is diagonally dominant (proposition (\ref{titi}) ).}
\item{ $\mathbf U$ is a bidiagonal and diagonally dominant matrix. $\mathbf Z=\mathbf U^{-1}$ has columns vector $\mathbf {z}_k, \ k=1,...,n$ such that $z_{k-i,k}> z_{k-i-1,k}$  with  {$ \small{0<i\leq k-1}$}. (proposition (\ref{final}))  }

\item{ Fixed  an upper  bandwidth $q$, the entries $z_{i,j}$ with $j>i+q$ of $\mathbf Z$  are considered negligible.}
\item {The preconditioner $\mathbf P= \tilde {\mathbf Z}  \tilde {\mathbf {Z^t}}$ is built as:
  \begin{equation}
\tilde z_{i,j}=\left\{
\begin{array}{l}
z_{i,j} \quad  {if \ \ j \leq i+q}\\[5mm]
0 \quad {if \ \ j > i+q}
\end{array}\right.
\end{equation}}
\item{ The  sparse factor $\tilde {\mathbf Z}$ is computed by $T$-orthogonalization procedure posing the sparsity pattern  S=$\{(i,j)\ / j>i+q \}$ }

}
\end{enumerate}
 $ \mathbf T$ is a diagonally dominant matrix then the incomplete inverse factorization of $\mathbf T$ exists for any choice of the sparsity pattern S on $\mathbf Z$~\cite{book4-16}.
\noindent  From computationally point of view, the $T$-orthogonalization procedure with  the sparsity pattern  S is based on  matrix-vector operations with computational cost of $5(q+1)$ floating point operations. Moreover, for each column vector $\tilde {\mathbf z}_k$ of $ \tilde {\mathbf{Z}}$ we work only on its  $q+1$ components $\tilde z_k[k-q], \tilde z_k[k-q+1],..., \tilde z_k[k]$ with consequently  global complexity of $5q(q+1) O(n)$.

\section{Practical Considerations}

\noindent In this section we give some practical details on the elliptic solver implementation with FSAI preconditioner, on a generic GPU architecture. The matrices $\mathbf A$, $\tilde {\mathbf Z}$ and $\tilde {\mathbf Z^T}$ are stored with the special storage format Compressed Sparse Row (CSR). The FSAI  is performed in serial on the CPU and its building requires a negligible time on total execution of the elliptic solver. We show the implementation of the Algorithm 1 outlines on the GPUs \cite{bell,book4-10,eurosam}. 

\begin{algorithm}[h!]
	\caption{FSAI-PCG solver}%
	\begin{small}
		\begin{algorithmic}[1]
			\STATE $k=0$; $\quad \mathbf x_0= D_{i,j}^0=2D_{i,j}^{t-1}$, the initial guess;\\[0.5 mm]
			\STATE $\mathbf r_0 = \mathbf b-\mathbf A \mathbf x_0$;\\[0.5 mm]
			\STATE ${\mathbf s_{0}}= \tilde {\mathbf Z}  \tilde {\mathbf {Z^t}} {\mathbf r_0} $, with ${\mathbf{ P}}= \tilde {\mathbf Z} \tilde {\mathbf Z^{t}}$ the  FSAI preconditioner;\\[0.5 mm]
			\STATE $\mathbf d_0 = \mathbf s_0$;\\[0.5 mm]
			\WHILE{$\big({ \|\mathbf r_{k}\|}/{\|\mathbf b\|} > \epsilon \ .and.\ k \leq n \big)$} 
				\STATE $\mathbf q_k=\mathbf A \mathbf d_k$; $\quad \alpha_k ={ (\mathbf{s}_k,\mathbf r_k)}/({\mathbf d_k, \mathbf q_k})$; $\quad \mathbf x_{k+1} = \mathbf x_{k} +\alpha_k \mathbf d_k$;\\[0.5 mm]
			\STATE ${\mathbf r_{k+1}} = {\mathbf r_k}-\alpha_k {\mathbf q_k}$; $\quad {\mathbf s_{k+1}}=\tilde {\mathbf Z} \tilde {\mathbf {Z^t}} {\mathbf r_{k+1}} $; $\quad \beta_k ={({ \mathbf s}_{k+1}, \mathbf r_{k+1})}/{({ \mathbf s}_{k},\mathbf r_{k})}$;\\[0.5 mm]
				\STATE $\mathbf d_{k+1}={\mathbf r}_{k+1}+\beta_{k} \mathbf d_{k}$; $\quad k=k+1$; \\[0.5 mm]
			\ENDWHILE
		\end{algorithmic}
	\end{small}
\label{PCG}
\end{algorithm}

\begin{figure}
			\centering

			\includegraphics[scale=.5]{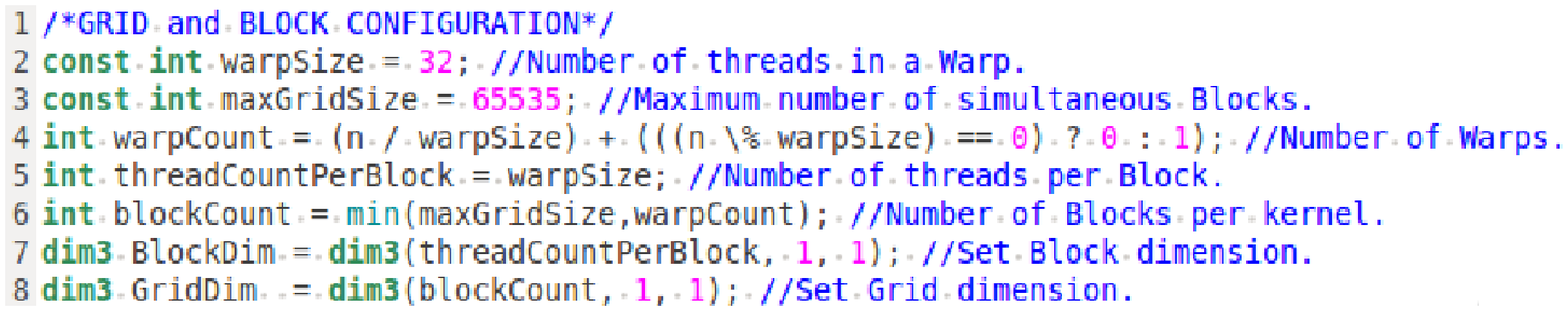}
			\caption{\small {Grid and block configuration.}}
			\label{fig:grid_configuration}
\end{figure}

\begin{figure}
			\centering
			\includegraphics[scale=.5]{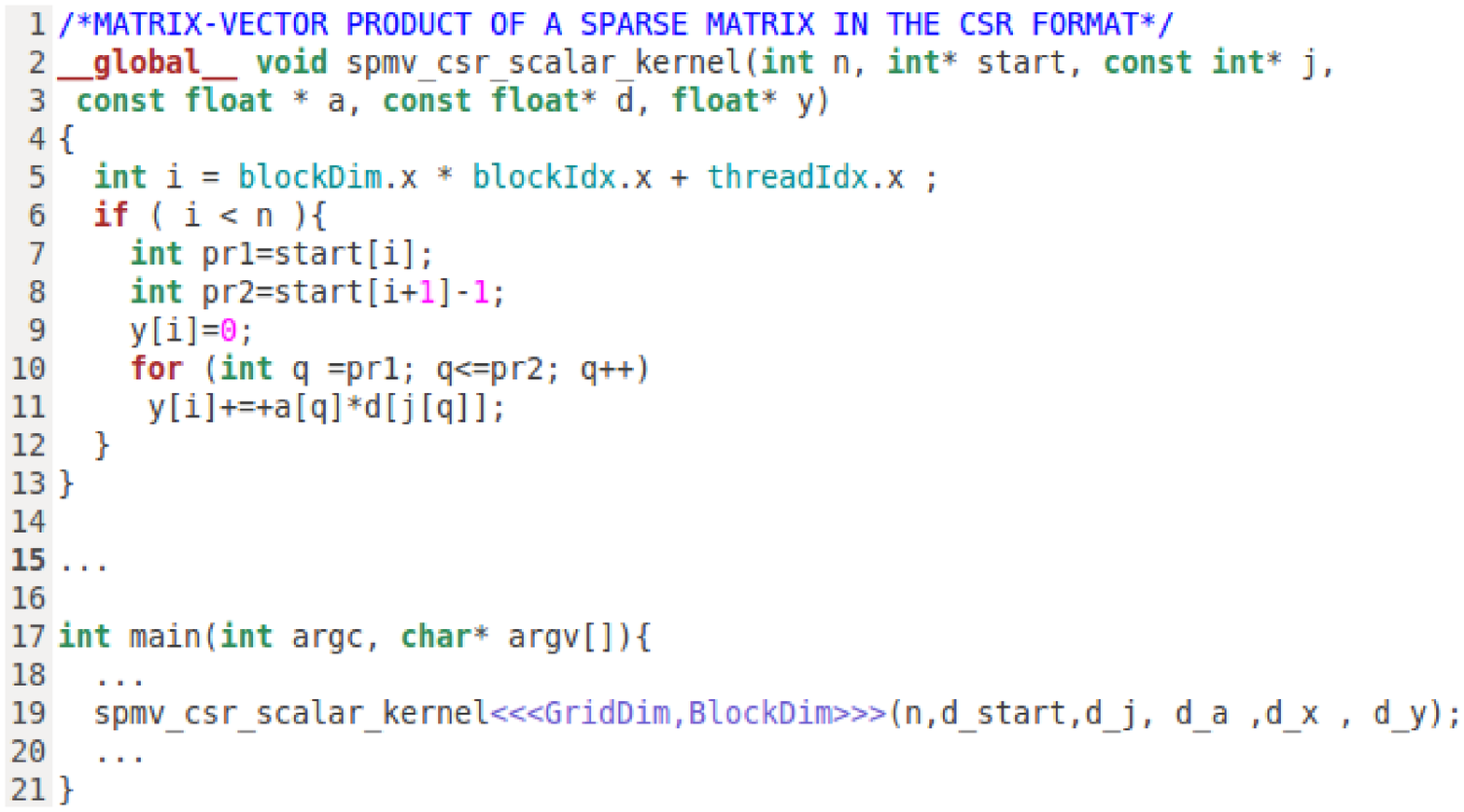}
			\caption{\small {Naive GPU implementation of a function for the matrix-vector product of a matrix in CSR format.}}
			\label{fig:spmv_csr_kernel}
\end{figure}

In details, our solver is implemented by means of the CUDA language with the auxiliary linear algebra libraries CUBLAS, for the ``dot product '' (\texttt{xDOT}), ``combined scalar multiplication plus vector addition'' (\texttt{xAXPY}), ``euclidean norm'' (\texttt{xNRM2}) and ``vector by a constant scaling'' (\texttt{xSCAL}) operations, and CUSPARSE for the sparse matrix-vector operations in the PCG solver.
\noindent The linear algebra scientific libraries are extremely helpful to easily implement a software on the GPU architecture. A ``by-hand'' implementation (see  Figure~\ref{fig:grid_configuration} and \ref{fig:spmv_csr_kernel}) of the solver without the library features in reported as a tedious GPU  programming example. In this type of coding, the manually configuration of the grid of thread blocks is necessary. For example if we use the TESLA S2050 the variables \texttt{warpSize} and \texttt{maxGridSize} (respectively at lines 2 and 3) have to be assigned. In details,\texttt{warpSize} indicates the number of threads (32) in a warp, which is a sub-division use in the hardware implementation to coalesce memory access and instruction dispatch; \texttt{maxGridSize} is the maximum number of simultaneous blocks (65535). Furthermore, \texttt{warpCount} (at line 4) represents the number of warps and it depends on the dimension \texttt{n} of the problem. In the end, variables \texttt{threadCountPerBlock} and \texttt{blockCount} (respectively at lines 5 and 6) are the parameter used for setting the grid and block configuration (lines 7 and 8). For example, if \texttt{n~=~10000} is the size of a vector, then \texttt{threadCountPerBlock~=~32} and \texttt{blockCount~=~313}.

\noindent The figure~\ref{fig:spmv_csr_kernel} shows the matrix-vector product, with the matrix stored in CSR format. We observe that this implementation requires a large amount of kernel functions, invoked by the ``host'' (CPU) and executed on the ``device'' (GPU). 

\noindent In the following, we will show how to implement the Algorithm 1 outlines by using library features. In order to use the CUBLAS library it is necessary inizialize it by means of the following instructions:

\begin{Verbatim}[frame=single]
cublasStatus stat;
cublasInit();
\end{Verbatim}

\noindent For the use of the CUSPARSE library two steps are necessary. The first one consists, as follow, in the library initialization:

\begin{Verbatim}[frame=single]
cusparseHandle_t handle=0;
cusparseCreate(&handle);
\end{Verbatim}

\noindent moreover, it is recalled the creation and setup of a matrix descriptor:

\begin{Verbatim}[frame=single]
cusparseMatDescr_t descra=0;
cusparseCreateMatDescr(&descra);
cusparseSetMatType(descra,CUSPARSE_MATRIX_TYPE_GENERAL);
cusparseSetMatIndexBase(descra,CUSPARSE_INDEX_BASE_ZERO);
\end{Verbatim}

\noindent The library  avoids to configure  the grid of the thread blocks and it allows to write codes in a very fast way. For example, at line 6 of the Algorithm~\ref{PCG} the computation of $\mathbf q_k= \mathbf {Ad}_k$ is required, and this operation can be made simply by calling the CUSPARSE routine \texttt{cusparseScsrmv()}, that performs the operation $\mathbf q = a \mathbf {A ∗ d} + b \mathbf q$ as follows:

\begin{Verbatim}[frame=single]
cusparseScsrmv(handle, CUSPARSE_OPERATION_NON_TRANSPOSE, n, n,
	  a, descra,  A,  start,  j,  d, b, q);
\end{Verbatim} 

\noindent In our context, \texttt{A}, \texttt{j} and \texttt{start} represent the  symmetric positive-definite matrix $\mathbf A$, stored in the CSR format. More precisely, the vector \texttt{A} denotes the non-zero elements of the matrix $\mathbf A$, \texttt{j} is the vector that stores the column indexes of the non-zero elements, the vector \texttt{start} denotes, for each row of the matrix, the address of the first non-zero element and \texttt{n} represents the row and columns  number  of the square matrix $\mathbf A$. The  constants $a$ and $b$ are assigned to $1.0$ and $0.0$ respectively. Moreover it happens that at line 6 of the Algorithm~\ref{PCG}, the computation of ${ (\mathbf{s}_k,\mathbf r_k)}$ is performed by means the CUBLAS routine for the dot product:
\begin{Verbatim}[frame=single]
alfa_num = cublasSdot(n, s, INCREMENT_S, r, INCREMENT_R);
\end{Verbatim}
The constants \texttt{INCREMENT$\_$S} and \texttt{INCREMENT$\_$R} are both assigned to $1$. Last operation of line 6 in Algorithm~\ref{PCG}, is the $\mathbf x_{k+1} = \mathbf x_{k} +\alpha_k \mathbf d_k$ for updating the solution and it is implemented by calling the CUBLAS routine for the saxpy operation:  
\begin{Verbatim}[frame=single]
cublasSaxpy(n, alfa, d, INCREMENT_D, x, INCREMENT_X);
\end{Verbatim}
The constants \texttt{INCREMENT$\_$D} and \texttt{INCREMENT$\_$X} are both assigned to $1$.
In addition, the computation of ${\mathbf{s_{k+1}}}=\tilde {\mathbf{Z}}  \tilde {\mathbf{Z^t}} {\mathbf{r_{k+1}}}$ at line 7  is the preconditioning step of the linear system~(\ref{sistema}) and it is computed  by means of two matrix-vector operations performed as:
\begin{Verbatim}[frame=single]
cusparseScsrmv(handle, CUSPARSE_OPERATION_NON_TRANSPOSE, n, n,
	a, descra, Z_t, start_Z_t, j_Z_t, r, b, zt);
cusparseScsrmv(handle, CUSPARSE_OPERATION_NON_TRANSPOSE, n, n,
	a, descra, Z_t, start_Z, j_Z, zt, b, z);
\end{Verbatim}
In details, in the first call of \texttt{cusparseScsrmv()}, ${\mathbf{zt}}=\tilde {\mathbf{Z^t}} {\mathbf{r_{k+1}}}$ and ${\mathbf{s_{k+1}}}=\tilde {\mathbf{Z}} {\mathbf{zt}}$ are computed. We have outlined just few of computational operations because the other  will be performed in the same way. The parameters \texttt{handle}, \scalebox{0.8}{\texttt{CUSPARSE\_OPERATION\_NON\_TRANSPOSE}} and \texttt{descra} are discussed in the NVIDIA report~\cite{book4-2bis} in more detailed way. In summary, we highlight how the use of the standard library, designed for the GPU architecture, allow to optimize the computational oceanographic  simulation model. 



\section{Numerical Experiments}
In this section we focus on the important numerical issues of our elliptic solver implemented with GPU architecture in single precision. The solver is tested on three grid size resolutions of the NEMO-OPA ocean model (Table~\ref{tab:gird}).

\begin{small}
\begin{center}
\begin{table}
\begin{center}
\scalebox{0.9}{
\begin{tabular}{ccc}
\hline
 Matrix Name & Size & Matrix non-zeros elements\\
\hline
ORCA-2 & $180\times149$ & $133800$ \\ \hline
ORCA-05 & $751\times510$ & $1837528$ \\ \hline 
ORCA-025 & $1442\times1021$ & $7359366$\\ 
\hline
\end{tabular}}
\end{center}
\caption{NEMO-OPA grid resolutions.}
\label{tab:gird}
\end{table}
\end{center}
\end{small}

\noindent In the Table~\ref{tab:a}, we compare the performance in terms of PCG iterations of the proposed  inverse bandwidth preconditioner $\mathbf P $ respect to $\mathbf P^{-1}$, that is the diagonal NEMO-OPA preconditioner. We fix an accuracy of $\epsilon = 10^{-6}$  on the relative residue $r=||\mathbf{Ax}-\mathbf{b}||/||\mathbf b||$ on the linear system solution. The experiments are carried out in the case of well-conditioned $\mathbf A$ matrix, corresponding to the geographical case of $\phi \approx 0$ and in the case of ill-conditioned $\mathbf A$  with $\phi \approx \pi/2$. We can observe as in the worst case with $n$ large and $\mathbf A$ ill-conditioned the PCG with $\mathbf P^{-1}$  has a very slow convergence with a huge number of iterations to reach the fixed accuracy. The experiment, reported in the Table~\ref{tab:a}, highlights the poor performance of the  $\mathbf P^{-1}$ for solving the Laplace elliptic problem~(\ref{sistema}) within OPA-NEMO. 

\begin{small}
\begin{center}
\begin{table}
\begin{center}
\scalebox{0.9}{
\begin{tabular}{ccccc}
\hline
 $\mathbf A$ Dimension 
  & $\qquad ${$\mathbf P$}\scalebox{0.7}{ $(\phi \approx 0)$ }
  &  $\qquad $$\mathbf P^{-1}$\scalebox{0.7}{$(\phi \approx 0)$ }
  & $\qquad ${$\mathbf P$}\scalebox{0.7}{ $(\phi \approx \pi/2)$}
  & $\qquad $$\mathbf P^{-1}$\scalebox{0.7}{ $(\phi \approx \pi/2)$}   \\
\hline
ORCA-2  & 271 & 460 & 8725   & 26820\\ \hline
ORCA-05 & 1128 & 1593 & 22447  & 86280\\ \hline 
ORCA-025 & 2458 & 3066 & 28513  & 139742\\ 
\hline
\end{tabular}}
\end{center}
\caption{Comparison between  $\mathbf P$ and $\mathbf P^{-1}$ in terms of Number of Iterations of the PCG in the case $\mathbf A$ is well-conditioned ($\phi \approx 0$) and $\mathbf A$  ill-conditioned ($\phi \approx \pi/2 $), varying the problem dimensions.}
\label{tab:a}
\end{table}
\end{center}
\end{small}

In the following, we show the performance in terms of PCG iterations of  $\mathbf P$ respect to the AINV Bridson Class preconditioners, that believe to CUSP library. To be more specific, let us consider the $\mathbf P_{B1}$ and $\mathbf P_{B2}$ Bridson's preconditioners, obtained by means of the $A$-orthogonalization method. The first is given by posing a (fixed) {drop tolerance} and by ignoring the elements below the fixed tolerance \cite{book4-21} and in the second one is predetermined the number of non-zeros elements on each its row.~\cite{book4-23}.

\begin{figure}[h!]
	\begin{minipage}[b]{0.5\linewidth}
		\begin{flushright}
			\includegraphics[scale=0.35]{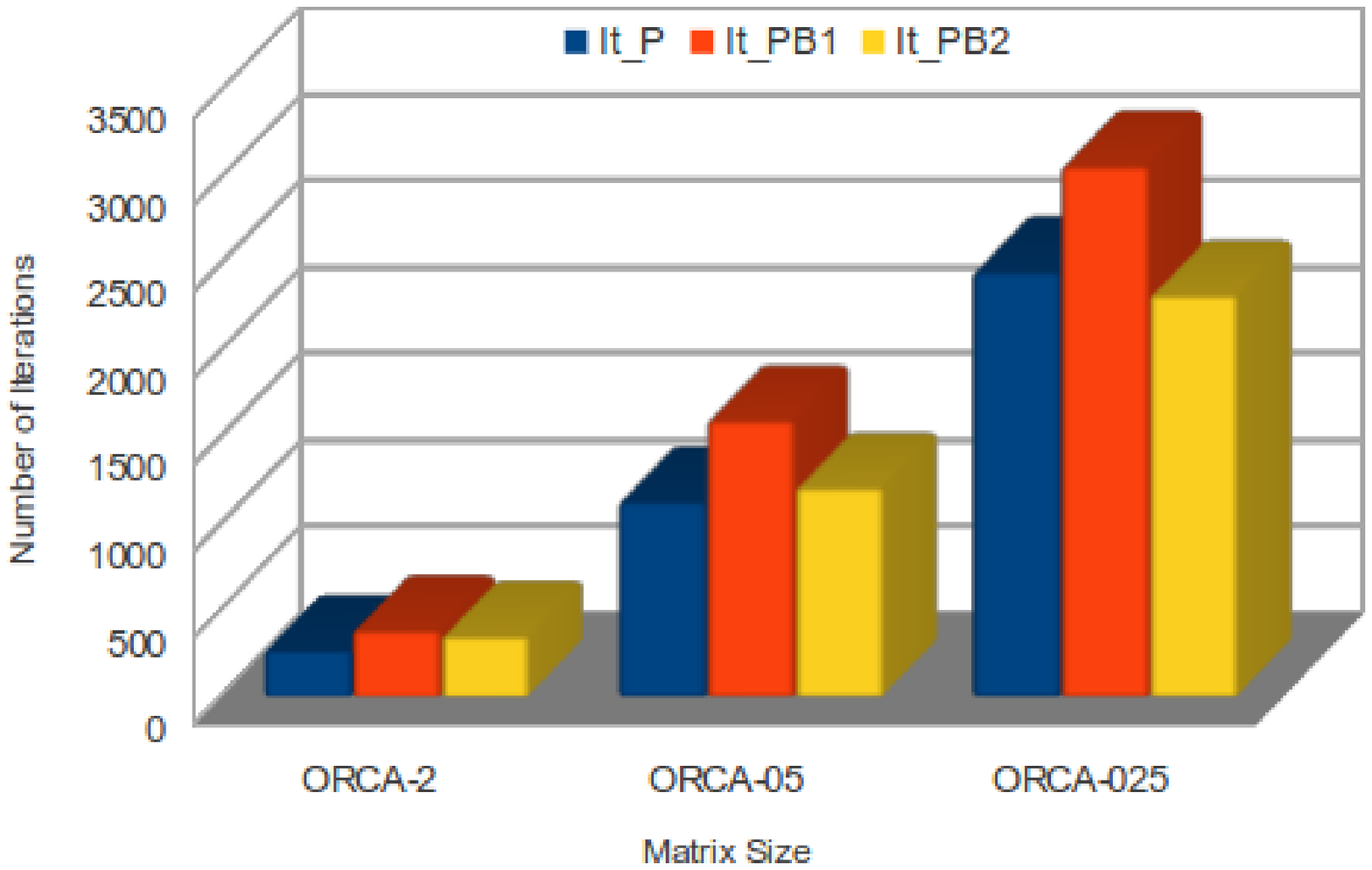}
			\caption{\small {Comparison between  $\mathbf P$, $\mathbf P_{B1}$ and $\mathbf P_{B2}$ in terms of Number of Iterations of the PCG ($y-$axis) when $\mathbf A$ is well-conditioned, varying the problem dimensions ($x-$axis) }}
			\label{fig:ITP_vs_Brids}
		\end{flushright}
	\end{minipage}
	\hspace{0.5cm}
	\begin{minipage}[b]{0.5\linewidth}
		\begin{flushleft}
			\includegraphics[scale=0.35]{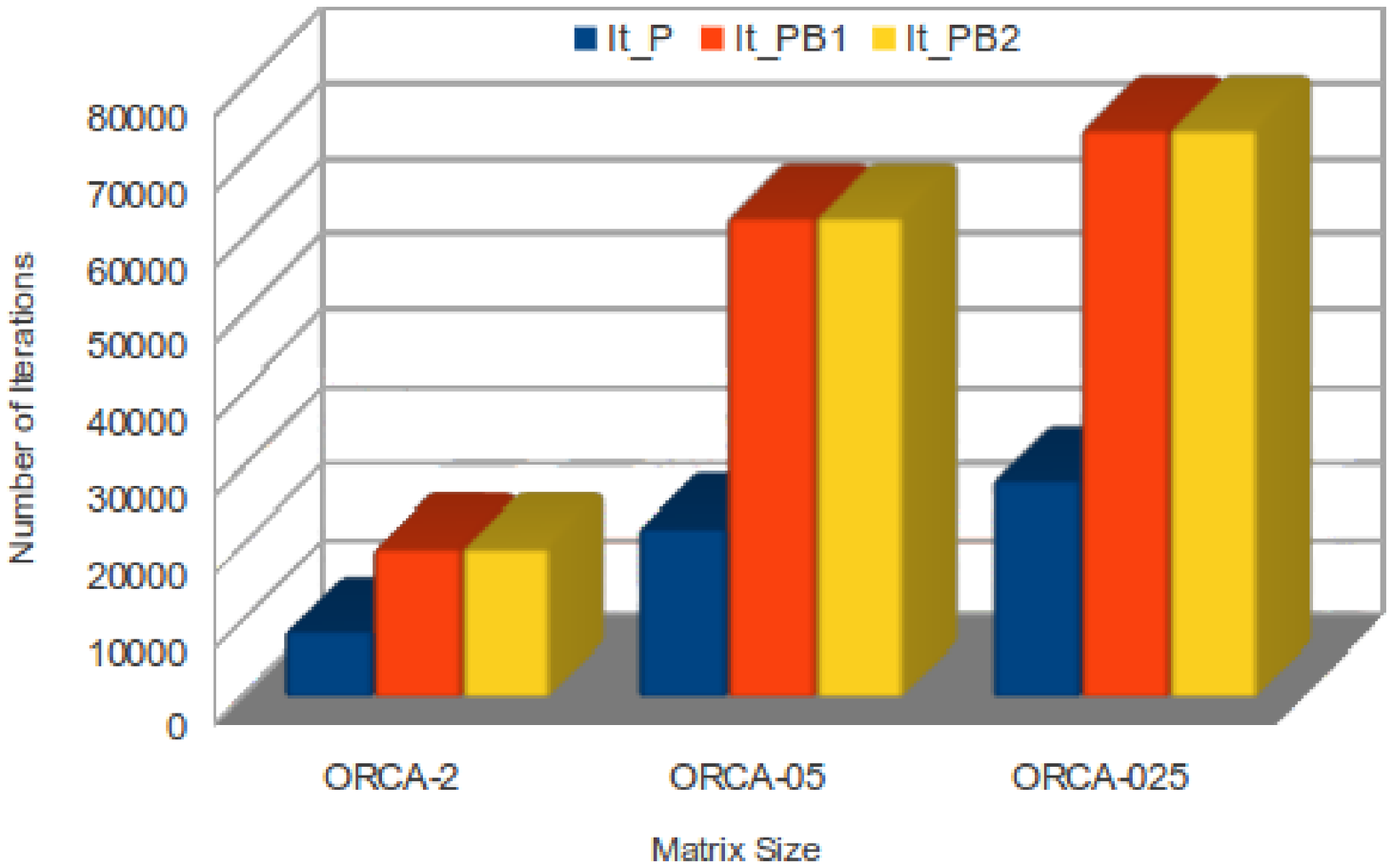}
			\caption{\small {Comparison between  $\mathbf P$, $\mathbf P_{B1}$ and $\mathbf P_{B2}$ in terms of  Number of Iterations of the PCG ($y-$axis) when $\mathbf A$ is ill-conditioned, varying the problem dimensions ($x-$axis) }}
			\label{fig:ITP_vs_Brids_geq}
		\end{flushleft}
	\end{minipage}
\end{figure}

\noindent The required accuracy on the solution is fixed to $\epsilon = 10^{-6}$ on the relative residue. In Figure~\ref{fig:ITP_vs_Brids} we report the PCG iterations of $\mathbf P$, $\mathbf P_{B1}$ and $\mathbf P_{B2}$ in the case of the matrix $\mathbf A$ well-conditioned $(\phi \approx 0)$. In Figure~\ref{fig:ITP_vs_Brids_geq} we present the case of the ill-conditioned $(\phi \approx \pi/2)$ matrix $\mathbf A$. \noindent The numerical results show as the number of iterations of the solver $\mathbf P$ is comparable to $\mathbf P_{B1}$ and $\mathbf P_{B2}$ when the dimensions of the problem are small or middle. Furthermore, it is strongly indicated to use $\mathbf P$ with a huge problem dimension. We test  $\mathbf {P}$, $\mathbf {P^{-1}}$,  $\mathbf P_{B1}$ and $\mathbf P_{B2}$ on the sparse matrix {\tt NOS6} of the Market Matrix database \cite{MaMa} by setting the required accuracy on the computed solution to $\epsilon = 10^{-6}$ and the band q of the preconditioner to 4. This sparse matrix is obtained in the Lanczos algorithm with partial re-orthogonalization Finite difference approximation to Poisson's equation in an L-shaped region, mixed boundary conditions.
\begin{figure}
	\centering
	\includegraphics[scale=0.39]{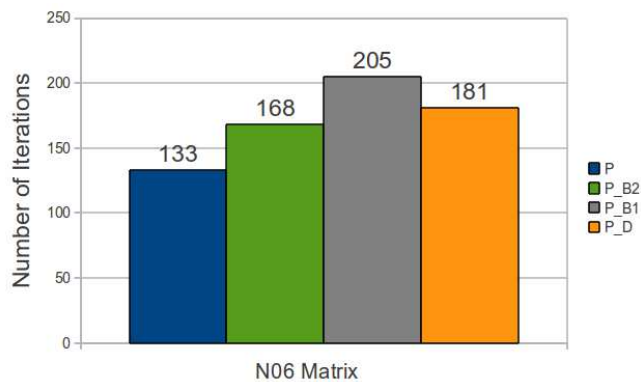}
	\caption{\small PCG iterations. P is the proposed preconditioner, P\_d is $\mathbf P^{-1}$, P\_B1  P\_B2 are $\mathbf P_{B1}$ and $\mathbf P_{B2}$ }
	\label{fig:N06}
\end{figure}

The figure~\ref{fig:N06} shows as $\mathbf {P}$ achieves the best performance in terms of iterations.\\
Finally, we test the elliptic solver implementation on GPU architecture. The numerical experiments are carried out on an ``NVIDIA TESLA S2050'' card, based on the ``FERMI GPU''. The ``TESLA S2050'' consists of 4 GPGPUs, each of which with 3GB of RAM memory and 448 processing cores working at 1.55 GHz. All runs are given on 1 GPU device. We have adopted CUDA release 4.0, provided by NVIDIA as a GPGPU environment  and the numerical code is implemented by using the single precision arithmetic. \\
As described in the previous sections by using scientific computing library it is not necessary  manually setting up the block and grid configuration on the memory device. The number of blocks required to store the elliptic solver input data (in CSR format) do not have to exceed the maximum sizes of each dimension of a GPU grid device. Schematic results of GPU memory utilization for ocean model resolutions are presented in the Table~\ref{tab:memory}. 
Observe that in our numerical experiments we do not fill the memory of the TESLA GPU and the simulations run also on cheaper or older boards, as for example the Quadro 4700FX. Generally, it is possible to grow the grid dimensions of the ocean model according to the memory capacity of the available GPU.

\begin{small}
\begin{center}
\begin{table}
\begin{center}
\scalebox{0.9}{
\begin{tabular}{ccc}
\hline
Matrix Name &  Non-zeros Elem. & Mem. Occ.\\
\hline
ORCA-2 &   $133800$ & $4$MB \\ \hline
ORCA-05 &  $1837528$ & $37$MB \\ \hline 
ORCA-025 & $7359366$ & $135$MB  \\ 
\hline
\end{tabular}}
\end{center}
\caption{Matrix memory occupancy. Mem Occ. is the full memory allocated memory on the GPU. }
\label{tab:memory}
\end{table}
\end{center}
\end{small}

\begin{figure}[h!]
	\begin{minipage}[b]{0.45\linewidth}
		\begin{flushleft}
			\includegraphics[scale=.6]{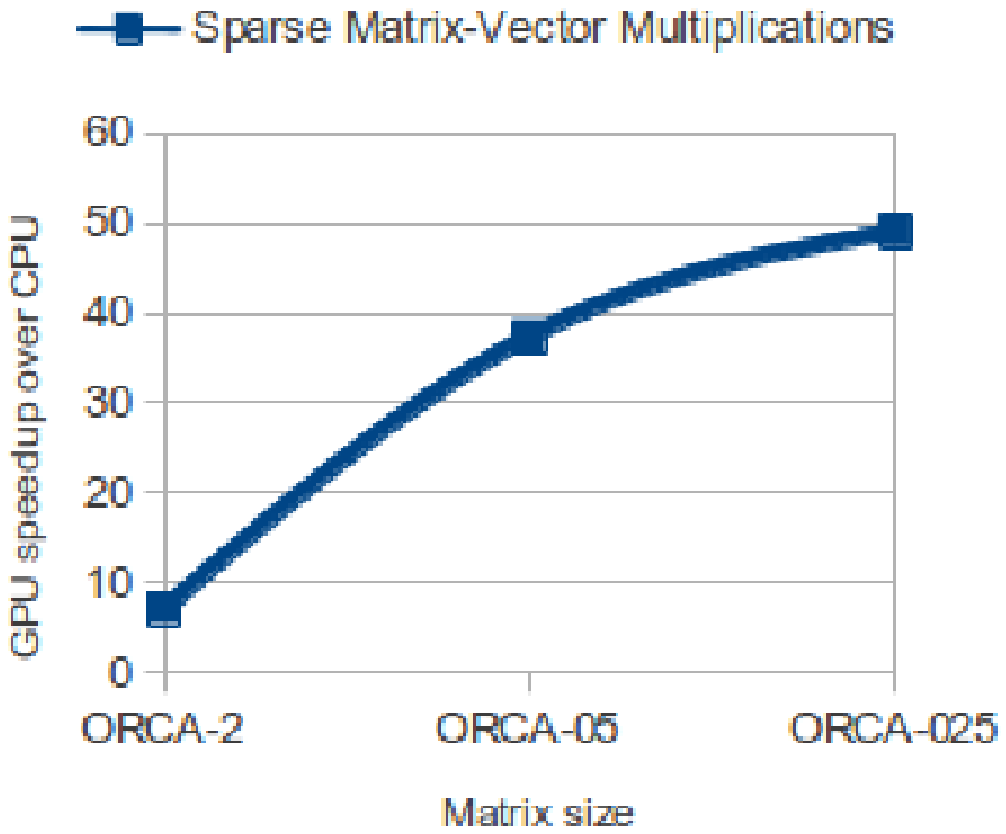}
			\caption{Sparse-Matrix Vector multiplications speed-up.}
			\label{fig:spmv}
		\end{flushleft}
	\end{minipage}
	\hspace{0.5cm}
	\begin{minipage}[b]{0.5\linewidth}
			\centering
			\includegraphics[scale=0.6]{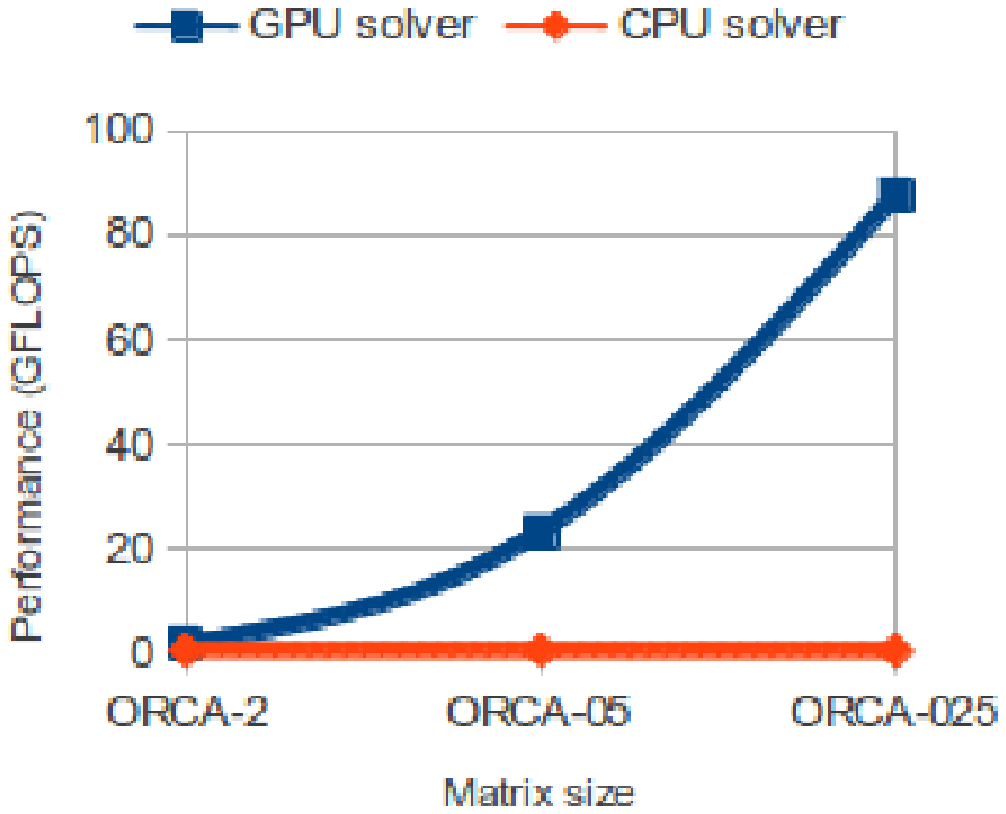}
			\caption{CPU and GPU comparison of the solver in terms of GFLOPS.}
			\label{fig:CPU_VS_GPU}
	\end{minipage}
\end{figure}

\noindent  The elliptic solver requires a large amount of Sparse-Matrix Vector (\texttt{cusparseCsrmv}) multiplications, vector reductions and other vector operation to be performed. CPU version is implemented in ANSI C executed in serial on  a 2.4GHz ``Intel Xeon E5620'' CPU, with 12MB of cache memory. Serial and GPU versions are in single precision. We test the performance of the solver in terms of Floating Point Operations (FLOPS). The performance of the numerical experiments (reported in the Figure~\ref{fig:CPU_VS_GPU}) are given in the case of $A$ ill-conditioned matrix ($\phi \approx \pi/2$). We count an average of the iterations of solver and the complexity of all linear algebra operations involved in both  serial and parallel implementations. The ''GPU solver`` (blue) and ''CPU solver`` (orange) curves represent the GFLOPS of the solver, respectively, for the CPU and GPU versions.

\noindent The main recalled computational kernels in the solver are the Sparse-Matrix Vector. From the Figure~\ref{fig:spmv} we highlight the improvement in terms of GFLOPs speed-up  by replacing \texttt{gemv()} with \texttt{cusparseCsrmv()} function. These results prove that, increasing the model grid resolution, it is possible to exploit the computational power of the GPUs. In details, the GPU solver implementation in the ORCA-025  configuration  has a peak performance of 87 GFLOPS respect to 1,43 GFLOPS of the CPU version.


\section{Conclusions}
The ocean modelling is a challenging application where expensive computational kernels are fundamental tools to investigate the physics of the ocean and the climate change. In a lot of applications, the elliptic Laplace equations are used in the complex mathematical models; they represents critical computational points since the  convergence of the numerical  solvers to a solution, within a reasonable number of iterations, it is  not always guaranteed. 
In our case, this happens to the preconditioning technique of the OPA-NEMO ocean model, for which we prove to be inefficient and inaccurate.
In this paper, we have proposed a new inverse preconditioner based on the FSAI method that shows better results respect to the OPA-NEMO diagonal one and to others of the Bridson class. Moreover, an important contribute is given by an innovative approach for parallelizing the elliptic solver on the Graphical Processing Units (GPU) by means of the scientific computing libraries. The library based implementation of the computing codes allows to optimize oceanic framework reducing the simulation times and to develop  computational solvers easy-to-implement.

\section{Acknowledgments}
The computing resources and the related technical support used for this work have been provided by CRESCO/ENEAGRID ``High Performance Computing infrastructure'' and its staff with  particular acknowledgments to the researcher  Marta Chinnici; CRESCO/ENEAGRID is funded by ENEA, the ``Italian National Agency for New Technologies, Energy and Sustainable Economic Development'' and by national and European research programs. See \url{www.cresco.enea.it} for more information.

\end{document}